\documentclass[12pt]{article}
\usepackage{epic,eepic}
\usepackage[dvips]{epsfig}
\usepackage{color}
\usepackage{amsmath, amscd, amssymb}
\usepackage{longtable,textcomp}
\usepackage{cite}
\usepackage{enumerate}
\usepackage{authblk}

\def\draw #1 by #2 (#3){
  \vbox to #2{
    \hrule width #1 height 0pt depth 0pt
    \vfill
    \special{picture #3}
    }
  }
\def\scaleddraw #1 by #2 (#3 scaled #4){{
  \dimen0=#1 \dimen1=#2
  \divide\dimen0 by 1000 \multiply\dimen0 by #4
  \divide\dimen1 by 1000 \multiply\dimen1 by #4
  \draw \dimen0 by \dimen1 (#3 scaled #4)}
  }

\begin{document}
\renewcommand{\labelenumi}{\theenumi}
\newcommand{\qed}{\mbox{\raisebox{0.7ex}{\fbox{}}}}
\newcommand{\ep}{\epsilon}
\newcommand{\la}{\lambda}
\newcommand{\G}{\Gamma}
\newcommand{\AL}{\mathbb{A}}
\newtheorem{theorem}{Theorem}
\newtheorem{example}{Example}
\newtheorem{remark}{Remark}
\newtheorem{problem}[theorem]{Problem}
\newtheorem{defin}[theorem]{Definition}
\newtheorem{definition}[theorem]{Definition}
\newtheorem{lemma}[theorem]{Lemma}
\newtheorem{corollary}[theorem]{Corollary}
\newtheorem{nt}{Note}
\newtheorem{proposition}[theorem]{Proposition}
\renewcommand{\thent}{}
\newenvironment{pf}{\medskip\noindent{\textbf{Proof}:  \hspace*{-.4cm}}\enspace}{\hfill \qed \medskip \newline}
\newenvironment{spf}{\medskip\noindent{\textbf{Sketch of proof}:  \hspace*{-.4cm}}\enspace}{\hfill \qed \medskip \newline}
\newenvironment{pft2}{\medskip\noindent{\textbf{Proof of Theorem~\ref{2}}:  \hspace*{-.4cm}}\enspace}{\hfill \qed \medskip \newline}
\newenvironment{defn}{\begin{defin}\em}{\end{defin}}{\vspace{-0.5cm}}
\newenvironment{lem}{\begin{lemma}\em}{\end{lemma}}{\vspace{-0.5cm}}
\newenvironment{cor}{\begin{corollary}\em}{\end{corollary}}{\vspace{-0.5cm}}
\newenvironment{thm}{\begin{theorem} \em}{\end{theorem}}{\vspace{-0.5cm}}
\newenvironment{pbm}{\begin{problem} \em}{\end{problem}}{\vspace{-0.5cm}}
\newenvironment{note}{\begin{nt} \em}{\end{nt}}{\vspace{-0.5cm}}
\newenvironment{exa}{\begin{example} \em}{\end{example}}{\vspace{-0.5cm}}
\newenvironment{rem}{\begin{remark} \em}{\end{remark}}{\vspace{-0.5cm}}
\newenvironment{pro}{\begin{proposition} \em}{\end{proposition}}{\vspace{-0.5cm}}

\def\diag{\text{diag}}
\def\al{\alpha}
\def\Ga{\Gamma}
\def\be{\beta}
\def\la{\lambda}
\def\ge{\geq}
\def\le{\leq}

\setlength{\unitlength}{12pt}
\newcommand{\comb}[2]{\mbox{$\left(\!\!\begin{array}{c}
            {#1} \\[-0.5ex] {#2} \end{array}\!\!\right)$}}
\renewcommand{\labelenumi}{(\theenumi)}
\renewcommand{\b}{\beta}
\newcounter{myfig}
\newcounter{mytab}
\def\mod{\hbox{\rm mod }}
\def\scaleddraw #1 by #2 (#3 scaled #4){{
  \dimen0=#1 \dimen1=#2
  \divide\dimen0 by 1000 \multiply\dimen0 by #4
  \divide\dimen1 by 1000 \multiply\dimen1 by #4
  \draw \dimen0 by \dimen1 (#3 scaled #4)}
  }
\newcommand{\Aut}{\mbox{\rm Aut}}
\newcommand{\w}{\omega}
\def\r{\rho}
\newcommand{\DbF}{D \times^{\phi} F}
\newcommand{\autF}{{\tiny\Aut{\scriptscriptstyle(\!F\!)}}}
\def\Cay{\mbox{\rm Cay}}
\def\a{\alpha}
\newcommand{\C}[1]{\mathcal #1}
\newcommand{\B}[1]{\mathbb #1}
\newcommand{\F}[1]{\mathfrak #1}
\title{Bounding the parameter $\beta$ of a distance-regular graph with classical parameters}

\author[a]{Chenhui Lv}
	\author[a,b]{Jack H. Koolen\footnote{J.H. Koolen is the corresponding author.}}
	\affil[a]{\footnotesize{School of Mathematical Sciences, University of Science and Technology of China, Hefei, 230026, People's Republic of China}}
	\affil[b]{\footnotesize{CAS Wu Wen-Tsun Key Laboratory of Mathematics, University of Science and Technology of China, Hefei, 230026, People's Republic of China}}

\date{\today}

\maketitle
\pagestyle{plain}
\newcommand\blfootnote[1]{%
		\begingroup
		\renewcommand\thefootnote{}\footnote{#1}%
		\addtocounter{footnote}{-1}%
		\endgroup}
	\blfootnote{E-mail addresses:  {\tt lch1994@mail.ustc.edu.cn} (C. Lv),  {\tt koolen@ustc.edu.cn} (J.H. Koolen).}

\begin{abstract}
Let $\Gamma$ be a distance-regular graph with classical parameters $(D, b, \alpha, \beta)$ satisfying $b\geq 2$ and $D\geq 3$.  Let $r=1+b+b^2+\cdots+b^{D-1}$.
In 1999, K. Metsch showed that   there exists a positive constant $C(\alpha,b)$ only depending on $\alpha$ and $b$, such that if  $\beta \geq C(\alpha, b)r^2$, then
either $\Gamma$ is a Grassmann graph or a bilinear forms graph.

In this work, we show that for $b\geq 2$ and $D\geq 3$, then there exists a constant $C_1(\alpha, b)$ only depending on $\alpha$ and $b$, such that if $\beta \geq  C_1(\alpha, b)r$, then either $\Gamma$ is a Grassmann graph, or a bilinear forms graph.
\\


\end{abstract}
\section{Introduction}
For undefined notions and more details see \cite{bcn89} and \cite{dkt}. We introduce distance-regular graphs with classical parameters in the next section. 
In this paper, we study distance-regular graphs with classical parameters $(D, b, \alpha, \beta)$. It is known that the parameter $b$ is an integer such that $b\neq 0,-1$ 
(Proposition~\ref{b_integer}).  C. Weng \cite{weng} classified the distance-regular graphs with classical parameters $(D,b,\alpha,\beta)$ and $b\leq-2$ under weak conditions, see Remark~\ref{weng} below. 
In this paper we consider the distance-regular graphs with classical parameters $(D,b,\alpha,\beta)$ and $b\geq 1$.  
The  distance-regular graphs with classical parameters $(D,b,\alpha,\beta)$ and $b =1$ were classified by Terwilliger, see \cite[Theorem 6.1.1]{bcn89}.  
In \cite{KLPY24}, Koolen et al. showed that for $b \geq 2$ and $D$ much larger than $b$, the parameter $\alpha$ is bounded by $b(b+1)^2 +1$. This was improved in \cite{kgly24}, where 
Koolen et al. showed that for $b \geq 2$ and $D \geq 9$, the parameter $\alpha$ is bounded by $b^2(b+1)+1$. Moreover, they showed that for $b=2$ and $D \geq 12$,
the parameter $\alpha$ is bounded by $2$. 
Let $r = 1 + b + \cdots + b^{D-1}$. 
In \cite[Corollary 1.3]{metsch}, Metsch showed that if $b \geq 2$ and $\beta \geq C(\alpha, b)r^2$, then the graph must be a bilinear forms graph or a Grassmann graph, where $C(\alpha, b)$ is a constant only depending on $\alpha$ and $b$.
Moreover, in the same paper,  he showed that, if $\alpha = b-1$, then 
$\beta < \frac{1}{2b-1}(2b^4 + 2b^3 + 2b^2 + b -1)r $ or the graph is a bilinear forms graph, and if $\alpha=b$, then $\beta < \frac{8}{3}(b^2+2b)r$ or the graph is a Grassmann graph. 

In this paper we show that if $b \geq 2$ and $\beta \geq C_1(\alpha, b)r$, where $C_1(\alpha, b)$ is a constant only depending on $\alpha$ and $b$, then the graph must be a bilinear forms graph, or a Grassmann graph.

Our main result is the following.
\begin{theorem}\label{main}
Let $\Gamma$ be a distance-regular graph with classical parameters $(D, b, \alpha, \beta)$, such that $b \geq 1$ and $D \geq 3$. Let $r = 1+b + b^2+ \cdots + b^{D-1}$. Then one of the following holds:
\begin{enumerate}
\item $\Gamma$ is a Johnson graph;
\item $\Gamma$ is a Hamming or a Doob graph;
\item $\Gamma$ is a halved cube;
\item $\Gamma$ is the Gosset graph with intersection array $\{27, 10, 1; 1,10, 27\}$;
\item $\Gamma$ is a Grassmann graph;
\item $\Gamma$ is a bilinear forms graph;
\item $b \geq 2$ and $\beta < \max\{2\frac{b+2}{2b+3}r(b+2)(\alpha b + b + \alpha), 2\frac{b+2}{2b+3}r((b+1)(b^2 + b +2)-3)\}.$
\end{enumerate}
\end{theorem}

\begin{remark}
\begin{enumerate}
\item For $\alpha \in \{b-1, b\}$, \cite[Corollary 1.3]{metsch} gives slightly better bounds. 
\item For $\alpha =b =2$, the bound of \cite[Corollary 1.3]{metsch} was slightly improved in \cite{KLG2024}. 
\item In \cite{GK2018+}, they showed that if a distance-regular graph with classical parameters $(D, b, \alpha, \beta)$ satisfies $\alpha = b \geq 2$, $D \geq 9$ and 
$\beta = b(b^{D-1} +b^{D-2} \ldots +1)$, then the graph has to be the Grassmann graph $J_b(2D, D)$. 
\item In \cite{GK2019}, they showed that if a distance-regular graph with classical parameters $(D, b, \alpha, \beta)$ satisfies $\alpha+1= b =2$, $D \geq 3$ and 
$\beta = b^D -1$, then the graph has to be the bilinear forms graph $H_2(D, D)$.
\item For $b \geq 2$, the twisted Grassmann graphs as found by Van Dam and Koolen in \cite{vDK} are not geometric and have $\alpha = b$ and $\beta = b^2r+b$. 
This shows that for 
$\beta = b^2r+b$ we can have Item (7) of Theorem~\ref{main}, and hence the bound in Item (7) is quite tight.
\end{enumerate}
\end{remark}

Our proof closely follows the ideas of Metsch as given in \cite{metsch}. We do not assume that $\alpha$ is an integer at most $b$. Note that for Metsch’s result, it is crucial that $\alpha$ is an integer at most $b$. We obtain this once we show that the distance-regular graph is 
geometric. The linear algebra idea in the proof of Proposition~\ref{alphaleqb} was first employed in 
\cite{LTK2021}. 

Outline of the paper: In the next section we give preliminaries and definitions. In Section 3, we discuss partial linear spaces.
In this section we give a sufficient condition for that a distance-regular graph is the point graph of a partial linear space, due to  Metsch. 
For the point graph of a partial linear space, we propose the PLS$(\gamma)$ and SPLS$(c,s)$ properties in Section 4, and give sufficient conditions for a distance-regular graph to 
satisfy these properties. Also, in this section, we give some structural results for a distance-regular 
graph satisfying these properties. In Section 5, we give sufficient conditions for a distance-regular graph to be geometric. 
In Section 6, we introduce the ELS property and establish the existence of geodetically closed strongly regular subgraphs.
In the last section we give the proof of the Theorem~\ref{main}.

\section{Definitions and preliminaries}
All the graphs considered in this paper are finite, undirected and simple. The reader is referred to \cite{bcn89, dkt} for more information. Let $\Gamma$ be a connected graph with vertex set $V(\Gamma)$. The {\em distance} $d(x,y)$ between two vertices $x,y\in V(\Gamma)$ is the length of a shortest path between $x$ and $y$ in $\Gamma$. The {\em diameter} $D=D(\Gamma)$ of $\Gamma$ is the maximum distance between
any two vertices of $\Gamma$. For each $x\in V(\Gamma)$, let $\Gamma_i(x)$ be the set of vertices in $\Gamma$ at distance
$i$ from $x$ ($0\leq i\leq D$). In addition, define $\Gamma_{-1}(x)=\Gamma_{D+1}(x)=\emptyset$. For the sake of simplicity, we denote $\Gamma_1(x)$ by $\Gamma(x)$.   For any vertex $x$ of $\Gamma$, the subgraph induced on $\Gamma(x)$ is called the {\em local graph} of $\Gamma$ at $x$, and we denote it by $\Delta_x$.  For a vertex $x$ of $\Gamma$, the  cardinality $|\Gamma(x)|$ of $\Gamma(x)$ is called the {\em valency} of $x$ in $\Gamma$. In particular, $\Gamma$ is {\em regular} with  valency $k$ if $k=|\Gamma(x)|$ holds for all $x\in V(\Gamma)$.  
We also define for vertices $x, y$ of $\Gamma$ at distance $i$ $(i= 1,2,\ldots, D)$ the set $C_i(x,y) $ by $C_i(x, y) := \Gamma_{i-1}(x) \cap \Gamma(y)$. For vertices $u, v$ of $\Gamma$ 
at distance $j$ $(j =0, 1, \ldots, D-1)$, we define the set $B_j(u, v) $ by $B_j(u, v) := \Gamma_{j+1}(u) \cap \Gamma(v)$.

For a set $T$ of vertices of a connected graph and a vertex $x$, we define the distance $d(x, T)$ by $d(x, T) = \min\{d(x, y) \mid y \in T\}$. 

The {\em adjacency matrix} $A=A(\Gamma)$ of $\Gamma$ is the matrix whose rows and columns are indexed by vertices of $\Gamma$ and the ($x, y$)-entry is $1$ whenever $x$ and $y$ are adjacent and $0$ otherwise. The {\em eigenvalues} of $\Gamma$ are the eigenvalues of its adjacency matrix $A$.

A {\em clique} of a graph $\Gamma$ is a set of mutually adjacent vertices of $\Gamma$. We sometimes also refer to a complete subgraph of $\Gamma$ as a clique. 
Let $K_{m,n}$ be the complete bipartite graph with $m$ and $n$ be the number of vertices in each partition.

A connected graph $\Gamma$ is called {\em edge-regular} with parameters $(n, k, a)$ if $\Gamma$ has $n$ vertices, is $k$-regular and any two adjacent vertices have exactly $a$ common
neighbours. A connected graph $\Gamma$ is called {\em amply-regular} with parameters $(n, k, a, c)$ if $\Gamma$ is edge-regular with parameters $(n, k, a)$ and any two vertices at 
distance 2 have exactly $c$ common neighbours.

\subsection{Distance-regular graphs}  A graph $\Gamma$ is called {\em{distance-regular}} if there exist integers $b_i, c_i$ $(0\leq i\leq D)$ such that for any two vertices 
$x, y \in V(\Gamma)$ with $d(x, y)=i$, there are precisely $c_i$ neighbors of $y$ in $\Gamma_{i-1}(x)$ and $b_i$ neighbors of $y$ in $\Gamma_{i+1}(x)$, where we define
$b_D=c_0=0$. Note that in this case $C_i(x,y)$ contains exactly $c_i$ vertices for $i =1,2,\ldots, D$ and $B_i(x, y)$ contains exactly $b_i$ vertices for $i =0,1,\ldots, D-1$.
In particular, any distance-regular graph  is regular with valency $k := b_0$.  We define $a_i := k-b_i-c_i$ for notational convenience.  
Note that $a_i=|\Gamma(y)\cap \Gamma_i(x)|$ holds for any two vertices $x, y$ with $d(x, y)=i$ $(0\leq i\leq D)$ and that the numbers $a_i, b_i$ and $c_i$ $(0\leq i\leq D)$ 
are called the {\em intersection numbers} of $\Gamma$.  For an eigenvalue $\theta$ of $\Gamma$, the sequence $(u_i)_{i=0,1,...,D}$ = $(u_i(\theta))_{i=0,1,...,D}$
satisfying $u_0$ = $u_0(\theta)$ = $1$, $u_1$ = $u_1(\theta)$ = $\theta/k$, and
\begin{align}\label{standardseq}
c_i u_{i-1} + a_i u_i + b_i u_{i+1} = \theta u_i ~(i=1,2,\ldots,D-1)
\end{align}
is called the {\em standard sequence} corresponding to the eigenvalue $\theta$ (see \cite[p.128]{bcn89}). 
The following lemma is well-known.

\begin{lemma}[cf.~{\cite[Proposition~4.4.1]{bcn89}}]\label{standard}
Let $\Gamma$ be a distance-regular graph, $\theta \neq k$ an eigenvalue of  $\Gamma$, say with multiplicity $m$ and let $(u_i)_i$ be the standard sequence corresponding to 
$\theta$. 
Then there exists a map $f: V(\Gamma) \rightarrow \mathbf{R}^m$, $x \mapsto \overline{x}$ such that the inner product
$(\overline{x}, \overline{y})= u_{d(x, y)}$. 
\end{lemma}

The map $f$ is called the {\em standard representation} of $\Gamma$ with respect to $\theta$. 

\subsection{Geometric distance-regular graphs}
The following lemma is called the Delsarte bound. We give a proof for it, as we need some consequences later in this paper.

\begin{lemma}[cf.~{\cite[Proposition~4.4.6]{bcn89}}]\label{Del}
Let $\Gamma$ be a distance-regular graph with diameter $D \geq 2$, valency $k$ and smallest eigenvalue $\theta_{\min}$, say with multiplicity $m$.
Let $(u_i)_i$ be the standard sequence corresponding to $\theta_{\min}$. 
Then the order $c$ of a clique $C$ in $\Gamma$ is bounded by $c \leq 1 + \frac{k}{-\theta_{\min}}$. 

Assume that $c = 1 + \frac{k}{-\theta_{\min}}$. 
Then any vertex of $\Gamma$ has distance at most $D-1$ to $C$. Let $x$ be a vertex with $d(x, C) = j$ for some $1\leq j \leq D-1$. Then the number of vertices $\phi_j(x)$ in $C$ at distance $j$ from $x$ satisfies:
\begin{equation}\label{decl}
(u_{j+1} - u_{j})\phi_j(x) = (1+ \frac{k}{-\theta_{\min}}) u_{j+1},
\end{equation}
and hence $\phi_j(x)$ does not depend on $x$, only on $j$.
\end{lemma}

\noindent
{\bf Proof.}
Let $f: V(\Gamma) \rightarrow \mathbf{R}^m, x \mapsto \overline{x}$ be the standard representation of $\Gamma$ with respect to $\theta_{\min}$. 
Let $\overline{C} = \sum \overline{u}$ where the sum is taken over all vertices $u$ of $C$.
Then the Gram matrix $Gr$ with respect to $\{\overline{u} \mid u$ is a vertex of $C\}$ is the matrix that satisfies
$$Gr = u_1 J + (1-u_1)I.$$ As $u_1 = \frac{\theta_{\min}}{k}$, and $Gr$ has row sum $1 + (c-1)u_1$ which is non-negative, as $Gr$ is positive semidefinite, we obtain
$c \leq 1 + \frac{k}{-\theta_{\min}}$, and equality holds if and only if $\overline{C} = \mathbf{0}$. Assume that equality holds. 
Let $x$ be a vertex of $\Gamma$ such that $d(x, C) = j$. Then the inner product $(\overline{x}, \overline{C})$ satisfies 
$$0 = (\overline{x}, \overline{C}) = u_j \phi_j(x) + u_{j+1} (c - \phi_j(x)), \mbox{ for } 1\leq j \leq D-1, \mbox{ and }$$
$$0 = (\overline{x}, \overline{C}) = u_D \phi_D(x), \mbox{ if } j=D.$$
As $u_D\neq 0$, we have $\phi_D(x) = 0$. This implies that there is no vertex of $\Gamma$ has distance $D$ to $C$. This shows the lemma. 
\qed

A clique of $\Gamma$ with $1+\frac{k}{-\theta_{\min}}$ vertices is called a {\em Delsarte clique}. 

A distance-regular graph $\Gamma$ is called \emph{geometric} if there exists a set of Delsarte cliques $\mathcal C$ in $\Gamma$ such that each edge lies in a unique $C \in \mathcal C$. In this case, we say that $\Gamma$ is a geometric distance-regular graph with respect to $\mathcal C$. Moreover, a partial linear space $X = (V(\Gamma), \mathcal C, \in)$ (see Section~\ref{plsdefinition} for the definition) is naturally induced, whose point graph is $\Gamma$, and hence the elements of $\mathcal C$ are also called lines.

Let $\Gamma$ be a geometric distance-regular graph with respect to $\mathcal C$. Let $C$ be a Delsarte clique and $x$ be a vertex at distance $j$ from $C$ for $0 \leq j \leq D-1$. 
Define $\phi_j$ be the number of vertices in $C$ at distance $j$ from $x$. 
For two vertices $x$ and $y$ at distance $j$ for $1 \leq j \leq D,$ let $\tau_j$ be the number of cliques in $\mathcal C$ containing $x$ at distance $j-1$ from $y$. 
The number $\tau_j$ $(j =1,2,, \ldots, D)$ does not depend on the pair $x, y$, only at their distance.

We can express the intersection numbers of $\Gamma$ in terms of the $\phi_j$’s and $\tau_j$’s as follows (see \cite[Lemma 4.1]{-m}):
\begin{equation}\label{eqgeom}
c_i = \tau_i \phi_{i-1} \ \ (i=1,2, \ldots, D), \ \ b_i= - (\theta_{\min} + \tau_i)\left(1+ \frac{k}{-\theta_{\min}} - \phi_i\right) \ \ (i = 1, 2, \ldots, D-1).
\end{equation}

We have the following properties for the numbers $(\phi_i)_i$ and $(\tau_i)_i$ for a geometric distance-regular graph.
\begin{lemma}[cf.~{\cite[Theorem 5.5]{BHK2007}}]\label{phiincrease}
Let $\Gamma$ be a geometric distance-regular graph with valency $k$ and diameter $D \geq 2$. Suppose $\phi_1 > 1$. Then 
\[
1 < \phi_1 < ... < \phi_{D-1}.
\]
In particular, $D \leq s$, where $s+1$ is the order of a Delsarte clique. 
\end{lemma}

\begin{lemma}[cf.~{\cite[Lemma 4.1]{bang2018}}]\label{tauincrease}

Fix an integer $m \geq 2$. Suppose that $\Gamma$ is a geometric distance-regular graph with diameter $D \geq 2$ and smallest eigenvalue $-m$. If $\phi_1 \geq 2$, then
\[
2 \leq \phi_1 \leq \tau_2 < ... < \tau_D = m.
\]
\end{lemma}

\begin{lemma}[cf.~{\cite[Lemma 5.2]{BHK2007}}]\label{2phi1}
Let $\Gamma$ be a distance-regular graph with valency $k$ and diameter $D \geq 2$, containing a Delsarte clique $C$ of size $s + 1$. Then the following hold:
\begin{enumerate}
    \item Let $i$ and $j$ be positive integers such that $i + j \leq D - 1$. Then $\phi_i + \phi_j \leq s + 1$.
    \item If $\phi_1 > \frac{s + 1}{2}$, then $D = 2$.
\end{enumerate}
\end{lemma}

\subsection{Distance-regular graphs with classical parameters}

A distance-regular graph $\Ga$ is said to have \emph{classical parameters} $(D,b,\al,\beta)$
if the diameter of $\Ga$ is $D$ and the intersection numbers of $\Ga$ can be expressed as follows:
\begin{equation}\label{bi}
 b_i=([D]-[i])
(\beta-\al [i]),~~  0\le i\le D-1,\nonumber
\end{equation}
\begin{equation}\label{ci}
c_i=[i](1+\al [i-1]),~~1\le i\le D,\nonumber
\end{equation}
where \begin{equation}\label{coeff}
[j]= \begin{cases}
    \frac{b^j-1}{b-1} & \text{if } b\neq1, \\
    \binom{j}{1} & \text{if } b=1. \\
  \end{cases}\nonumber
\end{equation}
From \cite[Corollary 8.4.2]{bcn89}, we know that the eigenvalues of $\Gamma$ are 
\begin{equation}\label{eigen}
[D-i](\beta- \alpha[i])-[i]=\frac{b_i}{b^i}-[i], ~~0\leq i\leq D.\nonumber
\end{equation}
We note that $c_2=(b+1)(\al+1)$ and that if $b\geq1$, then the eigenvalues $\theta_i=\frac{b_i}{b^i}-[i] (0\leq i\leq D)$ of $\Gamma$ are in the natural ordering, i.e., $k=\theta_0>\theta_1>\cdots>\theta_D$.

From the following result we know that the parameter $b$ of a distance-regular graph with classical parameters $(D,b,\alpha,\beta)$ is an integer.
\begin{proposition}[cf.~{\cite[Proposition~6.2.1]{bcn89}}]\label{b_integer}
Let $\Gamma$ be a distance-regular graph with classical parameters $(D,b,\alpha,\beta)$ and diameter $D\geq3$. Then $b$ is an integer such that $b\neq0,-1$. 
\end{proposition}

The following lemma is an easy consequence of Proposition~\ref{b_integer}.
\begin{lemma}\label{posalpha}
Let $\Gamma$ be a distance-regular graph with classical parameters $(D,b,\alpha,\beta)$ and diameter $D\geq3$, such that $b \geq 1$. Then $\alpha \geq 0$.
\end{lemma}
{\bf Proof.}
As $b$ is a positive integer and $c_2 = (\alpha +1)(b+1)$ is clearly a positive integer, we have $\alpha(b+1)$ is an integer.
Now the intersection number $c_3$ satisfies $c_3 = (1+b + b^2)(1+(1+b)\alpha) \geq 1$. This means $(1+b)\alpha +1 >0$, and as it is an integer, we see $(1+b)\alpha +1 \geq 1$.
This shows the lemma. \qed

\vspace{2mm}
The members of many of the infinite families of the known distance-regular graphs with unbounded diameter have classical parameters, like the Johnson graphs, Hamming graphs, Grassmann graphs and the bilinear forms graphs where the parameter $b$ is positive. 
\begin{remark}\label{weng}
\begin{enumerate}
\item For $b\leq-2$,  C. Weng \cite{weng} classified the distance-regular graphs with classical parameters $(D,b,\alpha,\beta)$ and $b\leq-2$. His result states: If $D\geq4$, $a_1\neq0$ and $c_2>1$, then one of the following holds: $(1)$ $\Gamma$ is the dual polar graph $^2A_{2D-1}(-b)$, $(2)$ $\Gamma$ is the Hermitian forms graph $Her_{-b}(D)$, $(3)$ $\alpha=\frac{b-1}{2}$, $\beta=-\frac{1+b^D}{2}$ and $-b$ is a power of an odd prime.
\item Tian et al. \cite{tian2024} give the current state of art on distance-regular graphs with classical parameters with $b \leq -2$.
\end{enumerate}
\end{remark}

Before we look at the Delsarte bound for distance-regular graphs with classical parameters, we will determine the standard sequence for the smallest eigenvalue of such a distance-regular graph.
Let $b\geq 2$ be an integer and let $[i] = \frac{b^i-1}{b-1}$ and $r = [D]$, as before.
\begin{lemma}\label{r}
Let $\Gamma$ be a distance-regular graph with classical parameters $(D, b, \alpha, \beta)$ with $b\geq 2$ and $D \geq 2$.
The standard sequence $(u_i)_i$ corresponding to the smallest eigenvalue $-r$ of $\Gamma$ satisfies:
$u_0 =1$ and $u_i = \frac{(-1)^i}{\beta}\frac{(1 +\alpha [1])(1 + \alpha[2]) \cdots (1+\alpha[i-1])}{(\beta-\alpha [1])(\beta- \alpha[2]) \cdots (\beta-\alpha[i-1])}$ for $i \geq 1$.
\end{lemma}

\noindent
{\bf Proof.}
We have $c_i = [i](1+\alpha[i-1]) $ for $1\leq i \leq D$, $a_i = [i](\beta-1+ \alpha(r - [i]-[i-1]))$ for $ 1 \leq i \leq D$ and $b_i = (r - [i])(\beta-\alpha[i])$ for $0 \leq i \leq D-1$. 
Note that $k = b_0 = \beta r$.
The $u_i$’s satisfy $u_0 = 1$, $u_1 = \frac{-r}{k} = \frac{-1}{\beta}$ and 
$c_iu_{i-1} + a_i u_i + b_i u_{i+1} = -r u_i$ for $1 \leq i \leq D-1$. The lemma follows by an easy induction. 
\qed

\vspace{2mm}
The following result is a direct consequence of Lemmas~\ref{Del} and \ref{r} for a distance-regular graph with classical parameters.
\begin{lemma}\label{clDel}
Let $\Gamma$ be a distance-regular graph with classical parameters $(D, b, \alpha, \beta)$ with $D \geq 3$ and $b \geq 2$.
Then the order $c$ of a clique $C$ in $\Gamma$ is bounded by $c \leq \beta +1$. If equality holds,
then the number $\phi_j$ for a vertex at distance $j$ from $C$ satisfies $\phi_j = 1 + \alpha[j]$ for $0\leq j \leq D-1$. In particular $\phi_1 = 1 + \alpha$. 
\end{lemma}

\noindent
{\bf Proof.}
As the smallest eigenvalue $\theta_{\min}$ of $\Gamma$ is equal to $-r$ where $r = [D]$, and $k = \beta r$ we obtain that any clique $C$ with order $c$ satisfies 
$c\leq 1+\frac{k}{r} = 1 + \beta$. If equality holds then by Equation~(\ref{decl}) we obtain $\phi_j(x) = 1 + \alpha[j]$, as the standard sequence corresponding to eigenvalue $-r$ satisfies 
$u_i = \frac{(-1)^i}{\beta}\frac{(1 +\alpha [1])(1 + \alpha[2]) \cdots (1+\alpha[i-1])}{(\beta-\alpha [1])(\beta- \alpha[2]) \cdots (\beta-\alpha[i-1])}$ for $i \geq 1$, by Lemma~\ref{r}.
This shows the lemma. \qed

\vspace{3mm}
Now we can calculate the $\phi_j$'s and $\tau_j$'s  of Equation~(\ref{eqgeom}) for a geometric distance-regular graph with classical parameters. 
\begin{lemma}\label{phi}
Let $\Gamma$ be a distance-regular graph with classical parameters $(D, b, \alpha, \beta)$ such that $b \geq 2$ and $D \geq 3$, and assume that 
$\Gamma$ is geometric. Then the numbers $\phi_j$'s and $\tau_j$'s of Equation~(\ref{eqgeom}) satisfy
\begin{equation*} \phi_j = 1 + \alpha[j], \ \ (j =0, 1, \ldots, D-1), \ \ \tau_j = [j] \ \ (j =1,2, \ldots, D), \end{equation*} where
$[j] = \frac{b^j -1}{b-1}$.
\end{lemma}
{\bf Proof.}
The expressions for the $\phi_j$'s  follow from Lemma~\ref{clDel} and for the $\tau_j$'s follow now from Equation~(\ref{eqgeom}). 
\qed

Now we first consider $\alpha =0$. 

\begin{proposition}\label{alpha=0}
Let $\Gamma$ be a distance-regular graph with classical parameters $(D, b, \alpha, \beta)$ such that $b \geq 1$, $\alpha =0$, $\beta>1$ and $D \geq 3$.
If $\Gamma$ is geometric, then $\Gamma$ is a dual polar graph or a Hamming graph. 
\end{proposition}
{\bf Proof.}
As $\alpha =0$ we have $a_1 = \beta-1 > 0$ and $\Gamma$ is locally the disjoint union of cliques of order $\beta$. By \cite[Theorem 9.4.4]{bcn89}, we are now done. 
\qed

\subsection{Strongly regular graphs}

A graph $\Gamma$ is called a {\em strongly regular graph} with parameters $(n, k, \lambda, \mu)$ if $\Gamma$ has $n\geq 2$ vertices, is $k$-regular and 
any two distinct vertices have 
$\lambda$ (resp. $\mu$) common neighbours depending on whether they are adjacent or not. For more information and for proofs of some the results we mention without proof see
\cite[Chapter 10]{GR2001} and \cite{SRG}. 

If $\Gamma$ or its complement is disconnected, we say $\Gamma$ is {\em imprimitive}, and {\em primitive} otherwise. Note that the only connected imprimitive non-complete 
strongly regular graphs are the complete multipartite graphs $K_{m \times s}$ with $m \geq 2, s\geq 2$. 
It is known that a connected non-complete $k$-regular graph is strongly regular  if and only if it has exactly three distinct eigenvalues $k > \sigma_1 > \sigma_2$. 
Note that a non-complete connected strongly regular graph with parameters $(n, k, \lambda, \mu)$ has non-integral eigenvalues only if $n = 4\mu+1, k = 2\mu,$ and $\lambda = \mu -1.$

We say that a connected non-complete strongly regular graph with parameters $(n, k, \lambda, \mu)$ has {\em classical parameters} $(b, \alpha, \beta)$ if $k = \beta(b+1)$, 
$\lambda = \beta-1+\alpha b$ and $\mu = (b+1)(\alpha +1)$. Note that a distance-regular graph with classical parameters $(2, b, \alpha, \beta)$ is just a strongly regular graph 
with classical parameters $(b, \alpha, \beta)$. Also note that any connected non-complete strongly regular graph with distinct eigenvalues $k > \sigma_1 > \sigma_2$ has classical 
parameters $(b, \alpha, \beta)$ where $b+1 \in \{-\sigma_1, -\sigma_2\}$.

We now reformulate Theorem 8.6.3 of \cite{SRG} in terms of terms of $(b, \alpha, \beta)$. Note that, if $b > 0$, then the smallest eigenvalue is equal to $-b -1$. We refer to \cite[Sections 8.4.2, 8.6]{SRG} for the definition of the block graph of a 
2-$(v, k,1)$-design and the Latin square graph $LS_{m}(n)$.

\begin{theorem}\label{BBN}
Let $b$ be an integer at least 2. Let $\Gamma$ be a primitive strongly regular graph with classical parameters $(b, \alpha, \beta)$. 
Let $f(b, \alpha) = \frac{1}{2}b(b+1)((\alpha+1)(b+1) +1) + \alpha$. 
Then $\beta \leq f(b, \alpha)$ or one of the following hold:
\begin{enumerate}
\item $\alpha = b$ and $\Gamma$ is the block graph of a 2-$(b\beta + b +1, b+1,1)$-design;
\item $\alpha = b-1$ and $\Gamma$ is a Latin square graph $LS_{b+1}(\beta +1)$.
\end{enumerate}
\end{theorem}

\section{Partial linear spaces}\label{plsdefinition}
An {\em incidence structure} is a tuple $(\mathcal{P}, \mathcal{L}, \mathcal{I})$ where $\mathcal{P}$ and $\mathcal{L}$ are non-empty disjoint sets and $\mathcal{I} \subseteq \mathcal{P} \times \mathcal{L}$. The elements of $\mathcal{P}$ and $\mathcal{L}$ are called {\em points} and {\em lines}, respectively. If $(p,\ell) \in \mathcal{I}$, then $p$ is said to be \emph{incident} with $\ell$. In this case, we also say that $\ell$ contains $p$ or that $p$ lies on $\ell$. The {\em order of a point} is the number of lines it is incident with and similarly for lines. The {\em point-line incidence matrix} of $(\mathcal{P}, \mathcal{L}, \mathcal{I})$ is the $|\mathcal{P}| \times |\mathcal{L}|$-matrix such that the $(p, \ell)$-entry is 1 if $p$ is incident with $\ell$ and 0 otherwise. If it is clear from the context what the incidence relation is, then we omit it.

A \emph{partial linear space} is an incidence structure  such that each pair of distinct points are both incident with at most one line.

Let $X=(\mathcal{P}, \mathcal{L}, \mathcal{I})$ be a partial linear space. For a non-incident pair $(x, \ell)$ of a point $x$ and line $\ell$, let $\tau(x, \ell) $ be the number of points on $\ell$ that are collinear with $x$. 
We define $\tau(X)$ for a partial linear space $X=(\mathcal{P}, \mathcal{L}, \mathcal{I})$, as the number $\tau(X) = \max\{ \tau(x, \ell) \mid x \in \mathcal{P}, \ell \in \mathcal{L}, \mbox{ and } (x,\ell)\notin \mathcal{I}\}$.

The \emph{point graph} $\Gamma$ of an incidence structure $(\mathcal{P}, \mathcal{L}, \mathcal{I})$ is the graph with vertex set $\mathcal{P}$ and two distinct points are adjacent if 
and only if they are on a common 
line. Given a line $\ell$, the vertex set $\{ x \in \mathcal{P} \mid (x, \ell) \in \mathcal{I} \}$ induces cliques in $\Gamma$.
If the point graph is connected, then for a vertex $x$ and a line $\ell$,  we define the distance $d(x, \ell)$ by 
$d(x, \ell)=\min\{d(x, y) \mid y \in \ell \}$. For two vertices $x, y$ at distance 2 in $\Gamma$, we denote by $[x, y]$  the set of  lines $\ell$ through $x$ such that $d(y, u) \leq 2$ for all $u \in \ell$. By Lemma~\ref{Del}, the set of lines $[x, y]$ are exactly the lines through $x$ at distance 1 from $y$ for a geometric distance-regular graph. 

First we give the so-called \emph{claw bound}. It was first proven by Metsch \cite[Lemma 1.1]{metsch2} and was rediscovered by  Godsil \cite[Lemma 2.3]{godsil} and Koolen and Park \cite[Lemma 2]{kp}.
\begin{lemma}\label{clawbound}
Let $\Gamma$ be an amply-regular graph with parameters $(n, k, \lambda, \mu)$.
Assume that there exists a positive integer $s$ such that the following holds:
 $$(s+1)(\lambda +1) - k > (\mu-1)\binom{s + 1}{2}.$$
 Then $\Gamma$ does not contain an induced $K_{1,s +1}$.
\end{lemma}

Now we recall the following result of Metsch. Note that the first condition is just the claw bound.

\begin{theorem}[cf.~{\cite[Result~2.1]{metsch}}]\label{metsch}
Let $\Gamma$ be an amply-regular graph with parameters $(n, k, \lambda, \mu)$.
Assume that there exists a positive integer $s$ such that the following two conditions are satisfied:
\begin{enumerate}
\item $(s+1)(\lambda +1) - k > (\mu-1)\binom{s + 1}{2};$
\item $\lambda+1  > (\mu -1)(2s -1)$.
\end{enumerate}
Define a line as a maximal clique with at least $\lambda+2 - (\mu -1)(s-1)$ vertices.
Then $X= (V(\Gamma), {\mathcal L}, \in)$ is a partial linear space, where $\mathcal L$ is the set of all lines, and $\Gamma$ is the point graph of $X$. Moreover, each vertex lies on at most $s$ lines.
\end{theorem}

%In \cite{KLG2024} they improved Theorem \ref{metsch} as follows. We give a proof of it for the convenience of the reader. 

\section{The SPLS$(s)$ Property}
We say that a distance-regular graph $\Gamma$ has the {\em PLS$(\gamma)$} property for some $\gamma\geq 3$ if $\Gamma$ is the point graph of a partial linear space 
$X= (V(\Gamma), {\mathcal L}, \in)$, where the set of lines $\mathcal L$ are the maximal cliques with at least $\gamma$ vertices.  

We have $\gamma\leq 1+\frac{k}{-\theta_{\min}}$ where 
$\theta_{\min}$ is the smallest eigenvalue of $\Gamma$, and equality if and only if $\Gamma$ is geometric.

Let $\Gamma$ be a distance-regular graph that is the point graph of a partial linear space $X= (V(\Gamma), {\mathcal L}, \in)$. A vertex $x$ of $\Gamma$ is called {\em Delsarte vertex} if all the lines in ${\mathcal L}$ through $x$ are Delsarte cliques. 

\begin{lemma}\label{pls}
Let $\Gamma$ be a distance-regular graph with classical parameters $(D, b, \alpha, \beta)$ such that $b \geq 2$ and $D \geq 3$.
Assume that $\Gamma$ satisfies the PLS$(\gamma)$ property for some $\gamma \geq 3$. Let $r= [D]$. 
Then a vertex $x$ lies on at least $r$ lines and it is on exactly $r$ lines if and only if it is a Delsarte vertex.

\end{lemma}
{\bf Proof.} By Lemma~\ref{clDel}, the order $c$ of a clique $C$ in $\Gamma$ satisfies $c \leq 1 + \beta$ with equality if and only if $C$ is a Delsarte clique. 
Now $k = \beta r$. By the definition of a partial linear space the lemma follows.
\qed

\vspace{2mm}
Now we prepare a lemma that will bound the parameter $\alpha$ of a distance-regular graph with classical parameters $(D, b, \alpha, \beta)$ such that $b \geq 2$.

\begin{lemma}\label{alpha}
Let $\Gamma$ be a distance-regular graph with classical prameters $(D, b, \alpha, \beta)$ such that $b \geq 2$ and $D \geq 3$.
Assume that $\Gamma$ satisfies the PLS$(\gamma)$ property for some $\gamma \geq 3$. 
Assume that $\Gamma$ contains two distinct Delsarte vertices at distance at most two. 
Then $\alpha \leq b$ and $\alpha$ is a non-negative integer.
\end{lemma}
{\bf Proof.}
Let $x, z$ be two vertices at distance two from each other such that $x$ is a Delsarte vertex and there exists a Delsarte clique $C'$ containing $z$ such that $x$ has at least one neighbour in $C'$ and hence $x$ has exactly $\alpha +1$ neighbours in $C'$. This means that $\alpha +1$ is a positive integer. Let $r= [D]$.  
Let $C_1, C_2, \ldots, C_r$ be the lines through $x$, all of which are Delsarte cliques.

The vertices $x$ and $z$ have $c_2 =(1+\alpha)(1+b)$ common neighbours and $z$ has exactly $\alpha +1$ or 0 neighbours in each of the $C_i$ for $i =1, 2, \ldots, r$.
This means that there are exactly $b+1$ $C_i$’s such that $z$ has $\alpha +1$ neighbours in $C_i$. Without loss of generality we may assume that $z$ has $\alpha +1$ neighbours in 
$C_i$ for $i=1, 2, \ldots, b+1$. On the other hand $x$ has exactly $\alpha +1$ neighbours in $C'$. Now the clique $C_i$ and $C'$ intersect in at most one vertex for $i =1, 2, \ldots, b+1$ 
and they do not intersect if $ i \geq b+2$. This means that $\alpha +1 \leq b+1$, and hence $\alpha \leq b$.

Let $y$ be a Delsarte vertex in $\Gamma$  at distance at most two from $x$ and distinct from $x$. 

Assume first that $x \sim y$. There exists a Delsarte clique $C'$ that contains $y$ but not $x$. As $C'$ is a maximal clique there exists a vertex $z$ in $C'$ not adjacent to $x$. 
So then we obtain $\alpha \leq b$ by the above.

Now assume $d(x,y) =2$. Then we take $z =y$, and we are done again. \qed

\vspace{2mm}
Now we will strengthen the PLS$(\gamma)$ property as follows:
\newline
We say that a distance-regular graph $\Gamma$ has the \emph{SPLS$(c,s)$ property} for some integers $c\geq 1$ and $s\geq 2$, if $\Gamma$ is the point graph of a partial linear space 
$X= (V(\Gamma), {\mathcal L}, \in)$, where the set of lines $\mathcal L$ are the maximal cliques with at least $a_1 +2 - (c-1)(s-1)$ vertices, $a_1  \geq (2s-1)(c-1)$,
each vertex lies on at most $s$ lines and $\tau(X) \leq c$. 

If $\Gamma$ has the SPLS$(c,s)$ property for some integers $c\geq 1$ and $s\geq 2$, then each line has at least $ a_1 +2 - (c-1)(s-1) \geq s(c-1)+2$ vertices, and clearly $c \leq c_2$ holds. If $c=c_2$, then we write SPLS$(s)$ instead of SPLS$(c_2, s)$.

Note that if a distance-regular graph $\Gamma$ satisfies the conditions of Theorem~\ref{metsch} for a positive integer $s$, then it has the SPLS$(s)$ property.
Note that the twisted Grassmann graphs $ \tilde{J}_q(2D+1, D)$ have the PLS($\frac{q^{D+1}-1}{q-1})$ property (see \cite{vDK}) , but do not have the SPLS$(s)$ property  for any $s\geq 2$.
This follows from the Proposition~\ref{alphaleqb}, as the twisted Grassmann graphs do not have Delsarte vertices. 

Now we give some bounds for $s$ for a distance-regular graph satisfying the SPLS$(s)$ property. 

\begin{lemma}\label{spls}
Let $\Gamma$ be a distance-regular graph with valency $k$ and smallest eigenvalue $\theta_{\min}$.
If $\Gamma$ has the SPLS$(s)$ property, then $-\theta_{\min} \leq s$, and if equality holds, then $\Gamma$ is a geometric distance-regular graph.
Let $\sigma$ be the smallest integer such that $\Gamma$ has the SPLS$(\sigma)$ property. Then 
$\sigma(a_1 +1) - k \leq(c_2-1)\binom{\sigma}{2}$, and $\sigma \geq -\theta_{\min}$.
\end{lemma}
{\bf Proof.}
As each clique has at most $1+\frac{k}{-\theta_{\min}}$ vertices, it is clear that $-\theta_{\min} \leq s$ and if equality holds, then each line is a Delsarte clique and hence the graph is geometric.

Let $x$ be a vertex that lies on exactly $\sigma$ lines. Then it is easy to show that $x$ lies in an induced $K_{1, \sigma}$, by an easy induction argument as any line has  at least $(c_2-1)\sigma+2$ vertices and that $\Gamma$ does not contain an induced $K_{1, \sigma +1}$. This shows the lemma, by applying Lemma~\ref{clawbound}. 
\qed

\vspace{2mm}
Now we give an extra condition on distance-regular graphs with classical parameters that have the SPLS$(s)$ property.

\begin{proposition}\label{alphaleqb}
Let $\Gamma$ be a distance-regular graph with classical parameters $(D, b, \alpha, \beta)$ such that $b \geq 2$ and $D \geq 3$.
Let $s$ be a positive integer and assume $\Gamma$ has the SPLS$(s)$ property. Let $r= [D]$.
Then the following hold:
\begin{enumerate}
\item $s \geq r$;
\item more than $\frac{c_2-2}{c_2-1} |V(\Gamma)|$ vertices are Delsarte vertices;
\item $\alpha$ is an integer such that $0 \leq \alpha \leq b$ and $\beta$ is a positive integer;
\item There exists two vertices $x, y$ at distance 2 such that $|[x, y]| \geq b+1$.
\end{enumerate}
\end{proposition}
{\bf Proof.}
The smallest eigenvalue of $\Gamma$ is equal to $-r$, and therefore $s \geq r$. 
By the SPLS$(s)$ property, let $X=(V(\Gamma), \mathcal{C}, \in)$ be the partial linear space with $\Gamma$ as its point graph,
where $\mathcal{C}$ is the set of maximal cliques that have at least $a_1 +2 -(s-1)(c_2 -1)\geq s(c_2-1)+2$ vertices. Note that all the Delsarte cliques of $\Gamma$ are in $\mathcal{C}$, as all cliques in $\Gamma$ have at most $\beta +1$ vertices and clearly 
$a_1 +2 -(s-1)(c_2 -1) \leq \beta+1$, as otherwise $\mathcal{C}$ would be empty.
Let $N$ be the point-line incidence matrix of $X$. 
Let $A$ be the adjacency matrix of $\Gamma$. 
Now $NN^T = A + Q$ where $Q$ is a diagonal matrix. Note that the rank of $N$ equals the rank of $A +Q$. 

By Lemma~\ref{pls} we see that a vertex $x$ of 
$\Gamma$ lies on exactly $r$ lines if and only if $x$ is a Delsarte vertex. This means that for all vertices $x$ of $\Gamma$ the entry $Q_{xx} \geq r$ and $Q_{xx}= r$ if and only if 
$x$ is a Delsarte vertex. 

Let $\delta$ be the dimension of the null-space of $A+Q$ and $\kappa$ be the number of Delsarte vertices. 
Let $\mathbf{u} \neq \mathbf{0}$ such that $(A+Q)\mathbf{u} = \mathbf{0}$. We find $\mathbf{u}^T(A+ rI)\mathbf{u} \geq 0$, as $A$ has smallest eigenvalue $-r$, and 
$\mathbf{u}^T(Q-rI)\mathbf{u} \geq \mathbf{0}$, as $Q_{xx} \geq r$ for all vertices $x$ of $\Gamma$. This implies that $A\mathbf{u} = -r \mathbf{u}$ and $Q\mathbf{u}= r\mathbf{u}$. 
Hence we obtain that the multiplicity $m$ of $r$ as an eigenvalue of $Q$ is at least $\delta$. Further we have that $m$ is exactly the number of vertices $x$ such that $Q_{xx}= r$.
This shows that $\delta \leq \kappa$. 

Now we determine the rank of $N$. This is at most the number of lines. So we need to estimate the number $\eta$ of lines. 
We do this by double counting the pairs $(x, \ell)$ where $x$ is a vertex of $\Gamma$ and $\ell$ is a line.
Let $w:=|\{(x, \ell) \mid x \in V(\Gamma), \ell \in \mathcal{C}, x \in \ell\}|.$
On the one hand we have that $w \leq s |V(\Gamma)|$ as every vertex lies on at most $s$ lines. 
On the other hand we find $w \geq \eta (a_1 +2 -(s-1)(c_2-1))$ as a line contains at least $a_1 +2 -(s-1)(c_2-1)\geq s(c_2-1)+2$ vertices.
Hence we obtain that
$$\eta \leq \frac{s|V(\Gamma)|}{a_1 +2 -(s-1)(c_2-1)}.$$
As $a_1 +2 -(s-1)(c_2-1) \geq s(c_2-1) +2 > s(c_2-1),$
we find that $\eta <\frac{1}{c_2-1}|V(\Gamma)|.$
So the rank of $N$ is less than $\frac{1}{c_2-1}|V(\Gamma)|$. 
This means that the dimension $\delta$ of the null-space of $A+Q$ is more than $|V(\Gamma)| -\frac{1}{c_2-1}|V(\Gamma)| = \frac{c_2-2}{c_2-1}|V(\Gamma)|$. 
As $b \geq 2$ and $\alpha \geq 0$, by Lemma~\ref{posalpha}, we have $c_2 = (\alpha+1)(b+1)\geq 3$. 
Then, as the number $\kappa$ of Delsarte vertices is at least the dimension $\delta$ and this number is in turn more than $|V(\Gamma)| -\frac{1}{c_2-1}|V(\Gamma)| = \frac{c_2-2}{c_2-1}|V(\Gamma)|$, we find more than half the vertices are Delsarte vertices. It means that there are adjacent vertices which are both Delsarte vertices. By Lemma~\ref{alpha} 
we find that $\alpha \geq 0$ is an integer and $\alpha \leq b$. 
If we take $x$ as a Delsarte vertex and $y$ a vertex at distance 2 of $x$, then $y$ has $\alpha+1$ or 0 neighbours on any line through $x$. As $c_2 = (\alpha +1)(b+1)$, we see that there are at least $b+1$ lines in $[x, y]$. This shows the proposition.
\qed

Now we give a sufficient condition for a distance-regular graph with classical parameters to have the SPLS$(s)$ property.

\begin{theorem}\label{alphabeta}
Let $\Gamma$ be a distance-regular graph with classical parameters $(D, b, \alpha, \beta)$ such that $b \geq 2$ and $D \geq 3$.
If $\beta \geq \frac{1}{3}(8\alpha b + 8b + 5 \alpha)r$, then $\Gamma$ has the SPLS$(\lfloor\frac{4r}{3}\rfloor)$ property. In particular, in this case, $\alpha$ is an integer such that $0 \leq \alpha \leq b$.
\end{theorem}
{\bf Proof.}
Assume that $\beta \geq \frac{1}{3}(8\alpha b + 8b + 5 \alpha)r$. We need to show that $\Gamma$ has the SPLS$(\lfloor\frac{4r}{3}\rfloor)$ property. By the definition of the SPLS$(s)$ 
property, we only need to show that the conditions in Theorem~\ref{metsch} are all satisfied for $s = \lfloor\frac{4r}{3} \rfloor$.
As $\lambda +1 = a_1 +1 = \beta+ \alpha(r-1)$, $\mu-1 = c_2 -1 = \alpha b + \alpha + b$ and $k = \beta r$, the conditions in Theorem~\ref{metsch} can be written as 
\begin{equation}\label{ab1}
\beta > \frac{s+1}{2(s+1-r)}\left((\alpha b + \alpha + b)s - 2\alpha(r-1)\right),
\end{equation}
and 
\begin{equation}\label{ab2}
\beta> (\alpha b + b + \alpha)2s - \alpha r - (\alpha+1)b.
\end{equation}
Let $s_0 := \frac{4r}{3}$ be a rational number. Then $$\frac{s+1}{s+1-r} < \frac{s_0}{s_0-r} = 4.$$

So to show that Conditions~(\ref{ab1}) and (\ref{ab2}) hold, it suffices to show that the following two conditions hold:
\begin{equation}\label{ab3}
\beta > 2(\alpha b + b + \alpha)s_0 - 4\alpha(r-1)= \frac{4}{3}(2(\alpha b + b + \alpha)r- 3\alpha(r-1))
\end{equation}
and
\begin{equation}\label{ab4}
\beta > 2(\alpha b + b + \alpha)s_0 - \alpha r - (\alpha +1)b= \frac{1}{3}(8(\alpha b + b + \alpha)r - 3\alpha r - 3(\alpha+1)b).
\end{equation}
When $\beta \geq \frac{1}{3}(8\alpha b + 8b + 5 \alpha)r$, Conditions~(\ref{ab3}) and (\ref{ab4}) clearly hold.
This shows the theorem. \qed
\section{Geometric distance-regular graphs}
In this section we will give a sufficient condition for a distance-regular graph with classical parameters $(D, b, \alpha, \beta)$ which has the SPLS$(s)$ property to be geometric. 

Our main result in this section is the following and we will give the proof using several claims. 
\begin{theorem}\label{geometric}
Let $\Gamma$ be a distance-regular graph with classical parameters $(D, b, \alpha, \beta)$ such that $b \geq 2$ and $D \geq 3$.
Let $s$ be a positive integer and assume $\Gamma$ has the SPLS$(s)$ property.
Let $r= [D]$. 
Assume that $\beta$ satisfies the following bounds:
\begin{enumerate}
\item $\beta > (b+2)(b+1)(s-1)-(b(b+1)+r-1)\alpha$;
\item $\beta > ((b+1)(b^2+b+1)-2)(s-b) -\alpha(b+1)r + \alpha(b+1)^2$;
\item $\beta >((b+1)(b^2 + b+1)-2+b) (s-b)  - \alpha(r-1) + b^2 -b $.
\end{enumerate}
Then $\Gamma$ is geometric. 
\end{theorem}
{\bf Proof.}
Note that by Proposition~\ref{alphaleqb}, we have $s \geq r$ and $0 \leq \alpha \leq b$ where $\alpha$ is an integer. Note that by (1), we have $\beta > (b^2+b+1)(b+1)$ as
$0\leq \alpha \leq b$ and $D \geq 3$. 

Now we will look at the line set $[x, y]$ for two vertices $x, y$ at distance 2. First, we show that this set has cardinality at most $2b+1$, and later, we will improve it to $b+1$.

\vspace{3mm}

\noindent
{\bf Claim 1}
If $\alpha \geq 1$, then$|[x, y]|< 2(b+1)$ for any two vertices $x$ and $y$ of $\Gamma$ at distance 2. If $\alpha = 0$, then$|[x, y]|\leq b+1 $ for any two vertices $x$ and $y$ of $\Gamma$ at distance 2.

\vspace{3mm}

\noindent
{\bf Proof of Claim 1.}

Let $M$  be the set of vertices consisting of neighbours of $x$ at distance at most 2 from $y$. Then 
\begin{align*}
|M|  & = k -b_2 = \beta r - \left(r - \frac{b^2-1}{b-1} \right) \left( \beta - \alpha\frac{b^2-1}{b-1} \right) \\ &= (b+1)(\beta + r\alpha- (b+1)\alpha) =(b+1)(a_1 +1 -b\alpha).
\end{align*}
As every line has at least $a_1 +2 -(s-1)(c_2-1)$ vertices, every line of $[x,y]$ has at least $a_1 +1 -(s-1)(c_2-1)$ vertices in $M$, because $x\not \in M$. 

Therefore 
\begin{align}\label{[x,y]}
|[x,y]| & \leq \frac{|M|}{a_1 +1 -(s-1)(c_2-1)} \\ & = \frac{(b+1)(a_1 +1 -b\alpha)}{a_1 +1 -(s-1)(c_2-1)} = (b+1)\left(1+\frac{(s-1)(c_2-1)-b\alpha}{a_1 +1 -(s-1)(c_2-1)}\right).
\end{align}

If $\alpha \geq 1$, by the SPLS(s) property we have $a_1 +1 \geq (2s-1)(c_2-1) +1= 2(s-1)(c_2 -1) +c_2$, and hence we have $|[x,y]| < 2(b+1)$. If $\alpha = 0$, then, by Condition (1) of the theorem, we have $$|[x, y]|\leq \frac{(b+1)\beta}{\beta-(s-1)b}<b+2.$$
This shows the claim.
\qed

\vspace{3mm}
\noindent
{\bf Claim 2}
\begin{enumerate}
\item We have $|[x, y]|\leq b+1$ for any two vertices $x$ and $y$ of $\Gamma$ at distance 2;
\item If a vertex has distance 1 to a line, then it has at most $b+1$ neighbours on the line;
\item $\Gamma$ has the SPLS$(b+1, s)$ property.
\end{enumerate}

\noindent
{\bf Proof of Claim 2.} 
By Claim 1, we may assume $\alpha \geq 1$. Let $q$ be the maximum integer such that there exist vertices $x, y$ at distance 2 with $|[x, y]| = q+1$. By Proposition~\ref{alphaleqb} we have $q \geq b$. 

We will first show that every line has at least $a_1 +2 - (s-1)q$ vertices. 
Let $\ell$ be a line and let $u\sim v$ be vertices on $\ell$. We'll show that any line $\ell_1$ through $u$ different from $\ell$ contains at most $q$ neighbours of $v$ other than $u$. 
Let $w$ on $\ell_1$ at distance 2 from $v$. Each neighbour of $v$ other than $u$ on the line $\ell_1$ gives a unique line in $[v, w]$. This means that $\ell_1$ contains at most $q$ neighbours of $v$ besides $u$. 
As $u$ lies on at most  $s$ lines, we obtain that $u$ and $v$ have at most $|\ell| -2 + (s-1)q$ common neighbours.
Hence $a_1 \leq |\ell|-2 + (s-1)q$, or equivalently, $|\ell| \geq a_1 +2 - (s-1)q$. 
Let $M$ be as in Claim 1. As every line has at least $a_1 +2 - (s-1)q$ vertices, every line of $[x, y]$ has at least $a_1 +1 -(s-1)q$ vertices in $M$. 
Therefore $(q+1)(a_1 +1 -(s-1)q) \leq |M|$. 
Let $f(x) = (x+1)(a_1 +1 -(s-1)x)$. Then $f(q) \leq |M|$. 
By Claim 1 we have $q \leq 2b$, and we have already seen $q \geq b$. 
We will show that $q =b$. Note that $|M| = (b+1)(a_1 +1 -b\alpha).$
The minimum of $f(x)$ for $b+1 \leq x \leq 2b$ is the minimum of $f(b+1)$ and $f(2b)$.

First we show that $f(b+1) \leq f(2b)$ holds.

Now $$f(2b) - f(b+1) = (2b+1)(a_1 +1 - (s-1)2b) - (b+2)(a_1 +1 - (s-1)(b+1)).$$ 
This implies that $$f(2b) - f(b+1) = (b-1)(a_1 +1) - (s-1)(3b^2 -b-2) = 
(b-1) (a_1 +1 -(s-1)(3b +2)).$$
So $f(2b) \geq f(b+1)$ if $a_1+1 \geq (s-1)(3b+2)$. 
We already found that $$a_1 +1 \geq 2(s-1)(c_2 -1) +c_2 = 2(s-1)((b+1)(\alpha +1)-1) + (b+1)(\alpha +1).$$
We see that $f(2b) \geq f(b+1)$ as $\alpha \geq 1$ and $b \geq 2$.
Now we need to compare $f(b+1)$ with $|M|$. 
We find 
\begin{align*}
f(b+1) - |M| & = (b+2)(a_1 +1 - (s-1)(b+1))- (b+1)(a_1 +1 -b\alpha) \\ & = a_1 +1 -(b+1)((b+2)(s-1) - b\alpha).
\end{align*}
As $a_1 + 1 = \beta + \alpha(r-1)$ and $\beta > (b+2)(b+1)(s-1)-(b(b+1)+r-1)\alpha$, we have
$$f(b+1) > |M|.$$
This means that $ q \leq b$ and hence this shows Item (1) of this claim.

Let $x$ be a vertex and $\ell$ a line such that $d(x, \ell)=1$. Then $x$ has at most $c_2$ neighbours on $\ell$. This means that there exists a vertex $y$ on $\ell$ at distance 2 from $x$.
There are at most $b+1$ lines in $[x, y]$ and each of them intersects in at most one vertex with $\ell$ and a line through $x$ that intersects with $\ell$ is in $[x, y]$. 
So $x$ has at most $b+1$ neighbours in $\ell$. 

Now it follows that  $\Gamma$ has the SPLS$(b+1, s)$ property, as through each vertex there at most $s$ lines. This shows the claim.
\qed

\vspace{2mm}
Now we will look at vertices at distance 3 in $\Gamma$.
\vspace{3mm}

 \noindent
 {\bf Claim 3}
 Let $x, y$ be two vertices at distance 3 in $\Gamma$, and $\ell$ a line through $y$. Then the number of vertices $z$ at distance 2 from $x$ such that $z \in \ell$ is less than $(b+1)(b^2 +b +1)$.  

\vspace{3mm}
\noindent
{\bf Proof of Claim 3.}
Let $U=\Gamma_2(x)\cap \ell$ and $W= \Gamma(x) \cap \Gamma_2(y)$. 
Then $|W| = c_3 = (b^2 +b +1)(1 + \alpha(b+1)$. 
If $u \in U$, then $\Gamma(u) \cap \Gamma(x) \subseteq W$. Hence every vertex of $U$ has $c_2$ neighbours in $W$.
Counting the number of edges between $U$ and $W$ we find  $|U|c_2 = \sum_{w \in W} |\Gamma(w) \cap \ell|.$
As a vertex outside $\ell$ has at most $b+1$ neighbours on $\ell$, we obtain $|U|c_2 \leq |W|(b+1)= (b+1)(b^2 + b +1)(1 + \alpha(b+1))$.
Now as $c_2 = (b+1) (\alpha +1)$, we see that $|U| \leq \frac{(b+1)(b^2 + b +1)(1 + \alpha(b+1))}{(b+1) (\alpha +1)} < (b^2 + b +1)(b+1)$, and this shows the claim.
\qed

\vspace{3mm}
\noindent
{\bf Claim 4}
If $\beta > ((b+1)(b^2+b+1)-2)(s-b) -\alpha(b+1)r + \alpha(b+1)^2$, then for any pair of vertices $x, y$ at distance 2, the number $|[x, y]|$ satisfies 
$|[x, y]| = b+1$.

\vspace{3mm}
\noindent
{\bf Proof of Claim 4.}
As above let $M$ be the set of neighbours of $x$ that have distance 1 or 2 from $y$. Then $|M| = k-b_2 = (b+1)\beta + \alpha(b+1)r - \alpha(b+1)^2$. 
Let $c := (b+1)(b^2+b+1).$
As $\Gamma$ has the SPLS$(b+1, s)$ propert, by Claim 2, we only need to show that $|[x, y]| \geq b+1$. Assume that $|[x, y]| \leq b$.
Any line through $x$ contains at most $\beta +1$ vertices, and therefore, there are at most $\beta$ vertices in $M$.

A line through $x$ that is not in $[x,y]$ contains a vertex at distance 3 from $y$, and hence by Claim 3 contains at most $c-2$ vertices in $M$.
As $x$ is on at most $s$ lines, it follows that
$|M| \leq |[x,y]|\beta + (s-|[x,y]|)(c-2) \leq b\beta + (s-b)(c-2)$, as $c-2 < \beta$. 
As $|M| = k-b_2 = (b+1)\beta + \alpha(b+1)r - \alpha(b+1)^2$, we obtain $\beta \leq  (s-b)(c-2) - \alpha (b+1)r + \alpha(b+1)^2$. 
So, if $\beta > (s-b)(c-2) - \alpha(b+1)r + \alpha(b+1)^2= ((b+1)(b^2+b+1)-2)(s-b) -\alpha(b+1)r + \alpha(b+1)^2$, then $|[x, y]| = b+1$, as we have 
$s \geq r \geq b^2 +b +1$ and $\alpha \leq b$. This shows Claim 4.
\qed
 
\vspace{3mm}
\noindent
{\bf Claim 5}
If $\beta > ((b+1)(b^2 + b+1)-2+b)(s-b)  - \alpha(r-1) + b^2 -b $, then for each pair of a vertex $x$ and a line $\ell$, with $d(x, \ell) =1$ and vertices $y_1, y_2$ on $\ell$ with $d(x, y_1) = d(x, y_2) =2$, we have that the sets of lines $[x, y_1]$ and $[x, y_2]$ 
are equal. 

\vspace{3mm}
\noindent
{\bf Proof of Claim 5.}
Assume that there are vertices $y_1, y_2 \in \Gamma_2(x)$ that are both on $\ell$ such that $[x, y_1]\neq[x, y_2]$. 
Put $S_i := \Gamma_i(x) \cap \ell$ and $n_i := |S_i|$ for $i=1,2.$
Further, put $T := [x, y_1]\cap[x, y_2]$ and $t :=  |T|$. Then $1 \leq n_1 \leq t \leq b$, as two distinct lines intersect in at most one vertex, $[x, y_1]\neq[x, y_2]$ and $|[x, y_1]|=|[x, y_2]| = b+1$.
Let $T_2$ be the set of lines through $x$ that are not in $T$. Then $t_2 := |T_2| \leq s-t$, since $x$ lies on at most $s$ lines. 
Let $\ell_2 \in T_2$ be a line through $x$ that is not in $[x, y_1]$. Then there exists a vertex $u$ on $\ell_2$ at distance 3 from $y_1$. Now we obtain $|\Gamma_2(u) \cap \ell| \leq c-1$ where $c = (b+1)(b^2+b+1)$. As $u \sim x$ and $d(u,y_1) =3$, we have $S_1 \subseteq \Gamma_2(u) \cap \ell$. As a consequence we find $|\Gamma_2(u) \cap S_2| \leq |\Gamma_2(u) \cap \ell| - |S_1| \leq c-1 - n_1.$
It follows that there are at most $c-1-n_1$ vertices $v \in S_2$ with $\ell_2 \in [x, v]$, as $d(u, v) \leq 2$ for such a vertex $v$. 

On the other hand, if $v \in S_2$, then $[x, v]$ contains at least $b+1-t$ lines of $T_2$. By counting pairs $(v, \ell_2)$ with $v \in S_2$ and $\ell_2 \in T_2 \cap [x, v]$ we obtain
$$n_2(b+1 -t) \leq |T_2|(c-1-n_1) \leq (s-t)(c-1-n_1).$$
As $t \leq b$ and $s \geq r > b+1$, we obtain $(b-t)(s-b-1) \geq 0$ which is equivalent with $s-t \leq (b+1 -t)(s-b)$. 
As a consequence we have $n_2 \leq (s-b)(c-1-n_1)$. This means that $|\ell| = n_1 + n_2 \leq n_1 + (s-b)(c-1-n_1) \leq 1 + (s-b)(c-2)$, as $n_1 \geq 1$ and $s \geq b+2$. 
By Claim 2, the graph $\Gamma$ has the SPLS$(b+1, s)$ property, and hence we have  $|\ell| \geq a_1 +2 -  (s-1)b$. 
Hence $\beta + \alpha(r-1) - (s-1)b+1 = a_1 +2 - (s-1)b \leq 1+ (s-b)(c-2)$ and therefore $\beta \leq (s-b)(c-2+b)  - \alpha(r-1) + b^2 -b $. 
So, if $\beta > (s-b)(c-2+b)  - \alpha(r-1) + b^2 -b $, the sets $[x, y_1]$ and $[x, y_2]$ are equal and this shows Claim 5.
\qed

\vspace{3mm}
\noindent
{\bf Claim 6}
Assume that $x$ is a Delsarte vertex and that $y$ is a neighbour of $x$. Then $y$ is also a Delsarte vertex.

\vspace{3mm}
\noindent
{\bf Proof of Claim 6.}
Let $C_1, C_2, \ldots, C_r$ be the lines through $x$ and without loss of generality, assume $y$ in $C_1$. The $C_i$’s are all Delsarte cliques by the fact that $x$ is a Delsarte vertex.
Let $C_1, D_2, \ldots, D_{\sigma}$ be the lines through $y$. Then $\sigma \geq r$, and if equality holds, then $y$ is a Delsarte vertex.
Let $\tilde{C}_i$ be the set of vertices in $C_i$ that are not neighbours of $y$, and $\tilde{D}_i$ be the set of vertices in $D_i$ that are not neighbours of $x$ for $i \geq 2$.
Let $\Delta$ be the bipartite graph with vertex set $\{\tilde{C}_2, \tilde{C}_3, \ldots, \tilde{C}_r, \tilde{D}_2, \tilde{D}_3, \ldots, \tilde{D}_{\sigma}\}$, and $\tilde{C}_i$ is adjacent to $\tilde{D}_j$
if and only if there are vertices $u \in \tilde{C}_i$ and $v \in \tilde{D}_j$ that are adjacent, for $i, j \geq 2$. 
We will obtain two bounds on the number of edges of $\Delta$ by double-counting.
Let $v_1, v_2$ be two distinct vertices of $\tilde{D}_j$ for some $j \geq 2$. By Claims 4 and 5, we have $[x, v_1] = [x, v_2]$ both with cardinality $b+1$. As ${C}_1 \in [x, v_1]$, we have that the valency $k_j$ of $\tilde{D}_j$ is at most $b$.
Let $C_i \in [x, v_1]$ for some $i \geq 2$. Then, as $C_i$ is a Delsarte clique, there are $\alpha +1$ neighbours of $v_1$ inside $C_i,$ by Lemma~\ref{clDel}. Note that $|\Gamma(y)\cap C_i-\{x\}|=\alpha$, therefore $\Gamma(v_1)\cap \tilde{C_i} \neq \emptyset$.
This means that $k_j$ is equal to $b$ and the number $\varepsilon$ of edges of $\Delta$ satisfies $\varepsilon = b(\sigma -1).$
Let $u_1, u_2$ be two distinct vertices of $\tilde{C}_i$ for some $i \geq 2$. Again by Claims 4 and 5, we have $[y, u_1] = [y, u_2]$ both with cardinality $b+1$. As ${C}_1 \in [y, u_1]$ we have that the valency $\kappa_i$ of $\tilde{C}_i$ is at most $b$. Therefore $\varepsilon \leq (r-1)b$ and hence $\sigma \leq r$. This means that $y$ is also a Delsarte vertex.
\qed

By Proposition~\ref{alphaleqb} we are done as we have at least one Delsarte vertex and hence all the vertices must be Delsarte vertices.
This shows that $\Gamma$ is geometric.
\qed

For geometric distance-regular graphs with classical parameters, we can provide a better bound for $\beta$ than the one in Claim 5 of Theorem~\ref{geometric}, using the same discussion as before. 

\begin{theorem}\label{geometric[]}
Let $\Gamma$ be a geometric distance-regular graph with respect to $\mathcal C$, having classical parameters $(D, b, \alpha, \beta)$, where $b \geq 2$ and $D \geq 3$. Let $r= [D]$. 
If $\beta > (r-\alpha -1)\alpha b +\alpha $, then for each pair of a vertex $x$ and a line $\ell$, with $d(x, \ell) =1$ and vertices $y_1, y_2$ on $\ell$ with $d(x, y_1) = d(x, y_2) =2$, we have that the sets of lines $[x, y_1]$ and $[x, y_2]$ are equal. 
\end{theorem}
{\bf Proof.}
Let $x,y$ be two vertices of $\Gamma$ at distance 2.
We will first show $|[x, y]| = b+1$. If $\ell$ is a line in $[x, y]$, then $d(y, \ell)=1$ as $\ell$ is a Delsarte clique,
and hence  $|[x, y]| = \tau_2 = b+1.$
Let $u, v$ be two vertices at distance 3 in $\Gamma$, and $m$ a line through $u$ and at distance 2 from $v$. Then the number of vertices $z$ at distance 2 from $v$ such that $z \in m$ is equal to $\phi_2= 
 1 + \alpha(b+1)$
and hence at most $b^2 +b +1$.

Let $\ell$ be a line at distance 1 from $x$.
Assume that there are vertices $y_1, y_2$ at distance 2 from $x$ that are both on $\ell$ such that $[x, y_1] \neq [x, y_2]$. Put, for $i=1,2$,  $S_i := \Gamma_i(x) \cap \ell$ and 
$n_i := |S_i|$.  By Lemma~\ref{clDel} we have $|\ell| = \beta +1$, $n_1 = \phi_1 = \alpha +1$ and $n_2 = |\ell| - n_1 = \beta -\alpha$, as $\Gamma$ is geometric.
This means that there are $\alpha +1$  lines through $x$ that intersect $\ell$, as two distinct lines intersect in at most one point and $n_1 = \alpha +1$. Those lines are in $T := [x, y_1] \cap [x, y_2]$. Let $t := |T|$. As $|[x, y_1]|= |[x, y_2]| = b+1$, and $[x, y_1] \neq [x, y_2]$, we have $\alpha +1 = n_1 \leq t \leq b$. 
As there are exactly $r$ lines through $x$, there are exactly $r -t$ lines through $x$ that are not in $T$. Let $u \sim x$ such that $d(y_1, u) =3$. Then the line $\ell_1$ through both $u$ and $x$ 
does not lie in $[x, y_1]$. Then $d(u, \ell) = 2$ and hence $|\Gamma_2(u) \cap \ell| = \phi_2 = 1+ \alpha(b+1)$, by Lemma~\ref{clDel}. 
As $u \sim x$ and $d(u, y_1) =3$, we find that every vertex in $S_1$ is at distance 2 from $u$. As a consequence 
$|\Gamma_2(u) \cap S_2| \leq |\Gamma_2(u) \cap \ell|  - |S_1| = 1+ \alpha(b+1) - (1+\alpha) = \alpha b$. 
It follows that there are at most $\alpha b$ vertices $v$ in $S_2$ such that $\ell_1 \in [x, v]$, as this implies $d(u, v) \leq 2$. 

On the other hand, for $v \in S_2$, the set $[x, v]$ contains at least $b+1 -t$ lines through $x$ that are not in $T$. 
By counting pairs $(v, \ell_2)$ with $v \in S_2$ and  lines $\ell_2$ in $[x, v]$ that are not in $T$ we have
$\beta- \alpha \leq (\beta-\alpha)(b+1-t) = n_2(b+1-t) \leq (r-t)\alpha b$, as $\ell_2 $ contains a vertex $w$ that is in $\Gamma_3(y_1) \cup \Gamma_3(y_2)$ and $t \leq b$. 

Therefore $\beta- \alpha \leq (r-t)\alpha b$. So, if $\beta> (r-\alpha -1)\alpha b +\alpha $, the sets $[x, y_1]$ and $[x, y_2]$ are equal. This shows the theorem. 
\qed

In the next theorem we give a bound on $\beta$ to insure that a distance-regular graph with classical parameters is geometric. 

\begin{theorem}\label{betabound}
Let $\Gamma$ be a distance-regular graph with classical parameters $(D, b, \alpha, \beta)$ such that $b \geq 2$  and $D \geq 3$.
Let $r = [D]$. If $$\beta \geq \max\{\frac{2b+4}{2b+3}r(b+2)(\alpha b + b + \alpha), \frac{2b+4}{2b+3}r((b+1)(b^2 + b +2)-3)\},$$ then $\Gamma$ is geometric.
\end{theorem}
{\bf Proof.}
We only need to show the existence a positive integer $s$ such that the conditions in both Theorem~\ref{metsch} and Theorem~\ref{geometric} are satisfied.
As we already have seen in the proof of Theorem~\ref{alphabeta}, Equations~(\ref{ab1}) and (\ref{ab2}), the conditions in Theorem~\ref{metsch} can be written as 
\begin{equation}\label{bb1}
\beta > \frac{s+1}{2(s+1-r)}\left((\alpha b + \alpha + b)s - 2\alpha(r-1)\right),
\end{equation}
and 
\begin{equation}\label{bb2}
\beta> (\alpha b + b + \alpha)2s - \alpha r - (\alpha+1)b.
\end{equation}
Let $s_0 := \frac{2b+4}{2b+3}r$ be a rational number and $s := \lfloor s_0 \rfloor$.
Then $\frac{s+1}{s+1-r} < \frac{s_0}{s_0 -r} = 2(b+2).$
So to show that Conditions~(\ref{bb1}) and (\ref{bb2}) hold, it suffices to show that the following two conditions hold:
\begin{align}\label{bb3}
\beta > (b+2)\left((\alpha b + b + \alpha)s_0 - 2\alpha(r-1)\right),
\end{align}
and
\begin{align}\label{bb4}
\beta &> 2(\alpha b + b + \alpha)s_0 - \alpha r - (\alpha +1)b.
\end{align}
But these two conditions are satisfied as $$\beta \geq \max\{\frac{2}{2b+3}r(b+2)^2(\alpha b + b + \alpha), 2\frac{b+2}{2b+3}r((b+1)(b^2 + b +2)-3)\}.$$
This shows that $\Gamma$ satisfies the SPLS$(s)$ property.
Now we need to consider the conditions in Theorem~\ref{geometric}.
These are 
\begin{equation}\beta > (b+2)(b+1)(s-1)-(b(b+1)+r-1)\alpha,\end{equation}\label{bb5}
\begin{equation}\beta > ((b+1)(b^2+b+1)-2)(s-b) -\alpha(b+1)r + \alpha(b+1)^2,\end{equation}\label{bb6}
and
\begin{equation}\beta > ((b+1)(b^2 + b+1)-2+b)(s-b)  - \alpha(r-1) + b^2 -b.\end{equation}\label{bb7}

As $$\beta \geq \max\{\frac{2}{2b+3}r(b+2)^2(\alpha b + b + \alpha), 2\frac{b+2}{2b+3}r((b+1)(b^2 + b +2)-3)\}$$ and $s_0 = \frac{2b+4}{2b+3}r$, we find that
\begin{equation}\beta > (b+2)(b+1)(s_0-1)-(b(b+1)+r-1)\alpha,\end{equation}\label{bb8}
\begin{equation}\beta > ((b+1)(b^2+b+1)-2)(s_0-b) -\alpha(b+1)r + \alpha(b+1)^2,\end{equation}\label{bb9}
and
\begin{equation}\beta > ((b+1)(b^2 + b+1)-2+b)(s_0-b)  - \alpha(r-1) + b^2 -b\end{equation}\label{bb10}
are satisfied and hence the conditions of Theorem~\ref{geometric} are satisfied. This shows that $\Gamma$ is geometric and this finishes the proof.
\qed

\vspace{2mm}
As a corollary we have the following consequence for distance-regular graph with classical parameters $(D, b, 0, \beta)$.

\begin{corollary}\label{alphazero} Let $\Gamma$ be a distance-regular graph with classical parameters $(D, b, \alpha, \beta)$, such that $b\geq 2$, $\alpha =0$ and $D \geq 3$. Let $r = [D]$. 
Then $\beta <  2\frac{b+2}{2b+3}(b^3 + 2b^2 + 3b -1)r$ holds.
\end{corollary}
{\bf Proof.}
For $\beta=1$ the corollary is clearly true, so we may assume $\beta >1$. If $\beta \geq 2\frac{b+2}{2b+3}(b^3 + 2b^2 + 3b -1)r$, then $\Gamma$ is geometric by Theorem~\ref{betabound}. By Proposition~\ref{alpha=0}, we find that $\Gamma$ is a dual polar graph or a Hamming graph. In the first case we have $\beta \leq b^2$, by \cite[Table 6.1]{bcn89}. In the last case we have $b=1$. 
 So the theorem follows.
\qed

\section{The ELS property and strongly regular subgraphs}
Let $\Gamma$ be a geometric distance-regular graph with respect to a set $\mathcal C$ of Delsarte cliques. 
We say that $\Gamma$ has the \emph{equal line set} (ELS) property if the following condition holds:
In the partial linear space $X = (V(\Gamma), \mathcal C, \in)$, for each pair $x \in V(\Gamma)$ and $\ell \in \mathcal C$ such that $d(x,\ell)=1$, and for every pair of distinct vertices $y_1$ and $y_2$ on $\ell$ at distance $2$ from $x$, the sets of lines $[x,y_1]$ and $[x,y_2]$ are equal.

%For each pair $(x, \ell)$, where $x$ is a vertex and $\ell$ a line with $d(x, \ell) =1$ and $y_1, y_2$ be two vertices on $\ell$ at distance 2 from $x$, the sets of lines $[x, y_1]$ and $[x, y_2]$ are equal. 

\begin{lemma}\label{vitallineset}
Let $\Gamma$ be a geometric distance-regular graph with respect to $\mathcal C$, having diameter $D \geq 3$.
Assume that $\Gamma$ satisfies the ELS property. Let $x,y$ be two vertices of $\Gamma$ at distance 2 in $\Gamma$.
Let $[x, y] = \{ \ell_1, \ell_2, \ldots, \ell_{\tau_2}\}$ and $[y, x] = \{m_1, m_2, \ldots, m_{\tau_2} \}$.
Let $z \in \Gamma(y) \cap \Gamma_2(x)$. Then $C_2(x, z) \subseteq \ell_1 \cup \ell_2 \cup \cdots \cup \ell_{\tau_2}$ if and only if there exists $i$ with 
$1 \leq i \leq \tau_2$ such that $z\in m_i$. 
\end{lemma}
{\bf Proof.} If $z \in m_i$, then by the ELS property $[x,y] = [x, z] = \{ \ell_1, \ell_2, \ldots, \ell_{\tau_2}\}$ and we are done.
Now assume that $C_2(x, z) \subseteq \ell_1 \cup \ell_2 \cup \cdots \cup \ell_{\tau_2}$. 
If $z$ is not on $m_1$, then there exists a vertex $u$ on $m_1$ adjacent to $x$, but not adjacent to $z$ as $z$ has exactly $\phi_1$ neighbours on $m_1$ including $y$, and  $x$ has 
has exactly $\phi_1$ neighbours on $m_1$. Hence there exists a vertex $w \in C_2(x, z)$ not adjacent to $y$ as both $y$ and $z$ have $c_2$ common neighbours with $x$.
Now by the ELS property we have $[y, x] = [y,w]$, and the line $L$ through $y$ and $z$ lies at distance 1 from $w$. So $L \in \{m_1, m_2, \ldots, m_{\tau_2}\}$. This shows that $z$ lies 
on $m_i$ for some $i =1, 2, \ldots, \tau_2$. This shows the lemma.
\qed

\begin{theorem}\label{srsg}
Let $\Gamma$ be a geometric distance-regular graph with respect to $\mathcal C$, having valency $k$, smallest eigenvalue $\theta_{\min}$, diameter $D \geq 3$, and $\phi_1 \geq 2$. Let $r = -\theta_{\min}$ and $\beta = \frac{k}{-\theta_{\min}}$.  
Assume that $\Gamma$ satisfies the ELS property. Then the following hold:
Let $x, y$ be a pair of vertices at distance 2 in $\Gamma$. Assume that $[x, y] = \{\ell_1, \ell_2, \ldots, \ell_{\tau_2}\}$. Let $\Sigma= \Sigma(x,y)$ be the subgraph induced on 
$\{x\} \cup \ell_1 \cup \ell_2 \cup \cdots \cup \ell_{\tau_2} \cup \{ z \in \Gamma_2 (x) \mid [x, y]= [x, z]\}$.
Then the graph $\Sigma$ has the following properties:
\begin{enumerate}
\item $\Sigma$ has diameter 2;
\item $\Sigma$ is geodetically closed;
\item If $\Sigma$ contains two vertices of a Delsarte clique $C \in {\mathcal C}$, then $\Sigma$ contains all the vertices of $C$;
\item $\Sigma$ is a geometric strongly regular graph with parameters 
\[
\left( \frac{\beta(\beta - \phi_1 + 1)(\tau_2 - 1)}{\phi_1} + \beta \tau_2 + 1, \ \beta \tau_2, \ \beta - 1 + (\phi_1 - 1)(\tau_2 - 1), \ \phi_1 \tau_2 \right);
\]
\item For every pair of vertices at distance 2, there exists a unique induced subgraph $\Delta$ of $\Gamma$ containing them such that $\Delta$ has diameter 2, is geodetically closed, and has the property that whenever $\Delta$ contains two vertices of a Delsarte clique in $\mathcal{C}$, it  contains the entire Delsarte clique. 
\end{enumerate}
\end{theorem}
{\bf Proof.} By Lemma~\ref{Del}, the number $\beta = \frac{k}{-\theta_{\min}}$ is a positive integer and each vertex lies on $r$ lines as $\Gamma$ is geometric.

Let $x,y$ be two vertices at distance 2 in $\Gamma$. Let $\ell_1, \ell_2, \ldots, \ell_r$ respectively $m_1, m_2, \ldots, m_r$ be the lines through $x$ respectively $y$. 
Assume  $[x,y] = \{\ell_1, \ell_2, \ldots, \ell_{\tau_2}\}$ and $[y,x] = \{m_1, m_2, \ldots, m_{\tau_2}\}$. Let $\Sigma$ be the graph as defined in the theorem. 
Then we have $y \in V(\Sigma)$ and $\{x\} \cup \Sigma(x) = \ell_1 \cup \ell_2 \cup \cdots \cup \ell_{\tau_2}$. 
By Lemma~\ref{vitallineset} we have that the neighbourhood of $z\in \Sigma_2(x)$ in $\Sigma$, denoted by $\Sigma(z)$ is the set $(\bigcup_{\ell \in [z,x]} \ell)\setminus \{z\}$. This, in 
particular, means that the valency of $z$ in $\Sigma$, denoted by $k_{\Sigma}(z)$ satisfies $k_{\Sigma}(z) = \beta \tau_2$.

\vspace{3mm}
\noindent
{\bf Claim 1}
Let $u$ be a vertex of  $\Sigma$, not adjacent to $y$.  Then $d_{\Sigma}(u, y) =2$ and $[y, u] = [y, x]$. 

\vspace{3mm}
\noindent
{\bf Proof of Claim 1.}
If $u \sim x$, then it is true by the ELS property.  So we may assume that $d(x, u) =2$. Assume that $[u, x] = \{p_1, p_2, \ldots, p_{\tau_2}\}$. By Lemma~\ref{vitallineset}, 
we have that the line $p_i$ lies completely inside $\Sigma$ for $i = 1, 2,\ldots, \tau_2$. If $y$ has a neighbour in each line $p_i$ for $i=1, 2, \ldots, \tau_2$ that is also adjacent to $x$, 
then it is evident that Claim 1 holds. Therefore, we may assume that there exists an $i =1, 2, \ldots, \tau_2$ such that all neighbours of $y$ in $p_i$ are not adjacent to $x$. Note that 
$x$ has $\phi_1 \geq 2$ neighbours in $p_i$ and they lie on $\phi_1$ lines of $[x,y] = \{\ell_1, \ell_2, \ldots, \ell_{\tau_2}\}$. 
Without loss of generality we may assume that there are vertices $w_j$ on $\ell_j$ for $j=1,2$ such that $\{w_1, w_2\} \subseteq p_i \cap \Sigma(x)$. 
By the ELS property, we have $[y, w_1]=[y, w_2] = [y, x]= \{ m_1, m_2, \ldots, m_{\tau_2}\}$. Let $[w_1, y] = \{q_1, q_2, \ldots, q_{\tau_2}\}$. Applying Lemma~\ref{vitallineset} to the tuple
$(y, w_1, \{m_t\}_{t=1}^{\tau_2}, \{q_t\}_{t=1}^{\tau_2})$ in the role of $(x, y, \{\ell_t\}_{t=1}^{\tau_2}, \{m_t\}_{t=1}^{\tau_2})$ shows that $w_2 $ is on a line $q_h$ for some 
$1 \leq h \leq \tau_2$. Note that the three vertices $u, w_1, w_2$ all lie on $p_i$ and that $w_1, w_2$ both lie on $q_h$. As each pair of vertices lie on at most one line we find that 
$p_i = q_h$ and therefore $u \in q_h$. This means that $d_{\Sigma}(u, y) = 2$. 
As $u$ and $w_1$ are both at distance 2 from $y$ and lie on the line $q_h$, we find $[y, u] = [y, w_1] = [y, x]$, by the ELS property. This shows Claim 1. 
\qed

\vspace{3mm}
\noindent
{\bf Claim 2}
Let $u \in \Sigma_1(x)$. Then there exists $w \in \Sigma_2(x)$ such that $d(w, u) = 2$. 

\vspace{3mm}
\noindent
{\bf Proof of Claim 2.}
Let $y \in \Sigma_2(x)$. Assume $y \sim u$, otherwise we are done. Consider the line $\ell$ through $x$ and $u$. Then $y$ has another neighbour $z$ on $\ell$ as $\phi_1 \geq 2$.
By Lemma~\ref{2phi1}, we have that $\phi_1 \leq \frac{\beta+1}{2}$, as $D \geq 3$. This means that the vertex $u$ has at most $\phi_1 -1\leq \frac{\beta-1}{2}$ neighbours  
on the line $m$ through $y$ and $z$ that are at distance two from $x$. But there are at least $\beta+1 - \phi_1 \geq \frac{\beta+1}{2}$ vertices on $m$ at distance 2 from $x$. 
So there exists a vertex $w$ on $m$ that is at distance 2 from both $u$ and $x$. This shows the claim. \qed

By the definition of $\Sigma(x, y)$ we have $\Sigma(x, y) = \Sigma(x, z)$ for all $z \in \Sigma_2(x)$. 

By Claim 1 we have $\Sigma(x, z) = \Sigma(z, x)$ for all $z \in \Sigma_2(x).$ Let $u \in \Sigma_1(x)$. Then there exists a vertex $w \in \Sigma_2(x)$ at distance 2 from $u$ by Claim 2. 
Hence $\Sigma(x, y) = \Sigma(x, w) = \Sigma(w, x) = \Sigma(w, u) = \Sigma(u, w)$ by Claims 1 and 2. 

It is now easily verified that $\Sigma$ is strongly regular and satisfies properties (1), (2), and (3).

For property (4), it is clear that $\Sigma$ has the given parameters. By direct calculation, the smallest eigenvalue of $\Sigma$ is $-\tau_2$, and hence a clique $C$ in $\Sigma$ is a Delsarte clique if and only if $|C|=\beta+1$. Therefore, each Delsarte clique $C \in \mathcal C$ with $C \subseteq V(\Sigma)$ is also a Delsarte clique in $\Sigma$. Together with property (3), this shows that $\Sigma$ is geometric.
 
 To show the last item, let $\Delta$ be an induced subgraph satisfying the conditions of Item (5) containing two vertices $x$ and $y$ at distance 2. 
 Then it follows easily that $\Delta$ is a strongly regular graph with the parameters as in Item (4). This means that for any vertex $z$ of $\Delta$ at distance $2$ of $x$ satisfies 
 $[x,z] = [x, y]$ and hence $\Delta = \Sigma(x,y)$.
 
This shows the theorem.
\qed

\begin{remark}
Note that the idea of finding the kind of subgraphs as in Theorem~\ref{srsg} goes back to Shult and Yanushka \cite{ShYa} and Brouwer and Wilbrink \cite{BrWi}. The proof of Theorem~\ref{srsg} is inspired by ideas from Gavrilyuk and Koolen~\cite{GK2019}.
\end{remark}

\section{Proof of Theorem~\ref{main}}
In this section we give the proof of Theorem~\ref{main}.

\vspace{2mm}
{\bf Proof of Theorem~\ref{main}.}
Terwilliger classified the distance-regular graphs with classical parameters $(D, b, \alpha, \beta)$ with $D \geq 3$ and $b=1$, and showed that they are the graphs in Items (1)-(4), see 
\cite[Theorem 6.1.1]{bcn89}.

Therefore we may assume that $b \geq 2$. 
If $\alpha=0$, then we are done by Corollary~\ref{alphazero}. So we assume $\alpha >0$.

We first consider the case that $\alpha \in \{b, b-1\}$. Assume that $\Gamma$ is not a Grassmann graph or a bilinear forms graph. In this case, by \cite[Corollary 1.3]{metsch}, we have: if $\alpha=b-1$, then 
\begin{align*}
\beta &< \frac{1}{2b-1}(2b^4 + 2b^3 + 2b^2 + b -1)r \\
&< \max\{\frac{2}{2b+3}r(b+2)^2(\alpha b + b + \alpha), 2\frac{b+2}{2b+3}r((b+1)(b^2 + b +2)-3)\};
\end{align*}

if $\alpha=b$, then
\begin{align*}
\beta &<\frac{8}{3}(b^2+2b)r \\
&<\max\{\frac{2b+4}{2b+3}r(b+2)(\alpha b + b + \alpha), \frac{2b+4}{2b+3}r((b+1)(b^2 + b +2)-3)\}.
\end{align*}

From now on, we may assume that $$\beta \geq\max\{\frac{2}{2b+3}r(b+2)^2(\alpha b + b + \alpha), 2\frac{b+2}{2b+3}r((b+1)(b^2 + b +2)-3)\}.$$

By Theorem~\ref{betabound}, the graph $\Gamma$ is geometric and hence it has Delsarte clique. This means that $1 \leq \alpha \leq b-2$ is an integer, as we already have considered
$\alpha \in \{0, b-1, b\}$.

Note that $\Gamma$ is geometric. Let $\mathcal C$ be a set of Delsarte cliques of $\Gamma$ such that $\Gamma$ is the point graph of the partial linear space $X=(V(\Gamma), {\mathcal C}, \in)$. 

By Theorem~\ref{geometric[]}, the graph $\Gamma$ has the ELS property. By Theorem~\ref{srsg}, there exists an induced subgraph $\Sigma$ of $\Gamma$ such that $\Sigma$ is a 
strongly regular graph with classical parameters $(b, \alpha, \beta)$.  Clearly $\Sigma$ is not a 
completely multipartite as the valency $k_{\Sigma} = (b+1)\beta > \mu(\Sigma) = (b+1)(\alpha +1)$ as the diameter of $\Gamma$ is at least 3. As $1 \leq \alpha \leq b-2$, we see that 
 $\mu(\Sigma) = (b+1)(\alpha +1) < (b+1)b$. Hence, by Theorem~\ref{BBN}, we have $\beta  \leq \frac{1}{2}(b+1)b((\alpha+1)(b+1)+1) + \alpha$. This gives a contradiction. This 
 shows the theorem. 
\qed

\bigskip
\begin{flushleft}

{\Large\textbf{Acknowledgments}}\vspace{0.2cm}

J.H. Koolen is partially supported by the National Natural Science Foundation of China (No. 12471335), and the 
Anhui Initiative in Quantum Information Technologies (No. AHY150000). We thank Bart De Bruyn and Ferdinand Ihringer for fruitful discussions.\

\end{flushleft}

\clearpage


\bigskip

\begin{thebibliography}{99}


\bibitem{bang2018} S. Bang, Diameter bounds for geometric distance-regular graphs, Discrete Math. 341 (2018), 253--260.

\bibitem{BHK2007} S. Bang, A. Hiraki, and J. H. Koolen, Delsarte clique graphs, European J. Combin. 28 (2007), 501--516.

\bibitem{bcn89} A. E. Brouwer, A. M. Cohen, and A. Neumaier, Distance-regular graphs, Springer-Verlag, Berlin, 1989.

\bibitem{SRG} A. E. Brouwer and H. Van Maldeghem, Strongly Regular Graphs, Cambridge University Press, Cambridge, 2022.

\bibitem{BrWi} A. E. Brouwer and H. A. Wilbrink, The structure of near polygons with quads, Geom. Dedicata 14 (1983), 145--176.

\bibitem{vDK} E. R. van Dam and J. H. Koolen, A new family of distance-regular graphs with unbounded diameter, Invent. Math. 162 (2005), 189--193.

\bibitem{dkt} E. R. van Dam, J. H. Koolen, and H. Tanaka, Distance-regular graphs, Electron. J. Combin. (2016) $\sharp$DS22.

\bibitem{GK2019} A. L. Gavrilyuk and J. H. Koolen, A characterization of the graphs of bilinear $(d \times d)$-forms over GF(2), Combinatorica 39 (2019), 289--321.

\bibitem{GK2018+} A. L. Gavrilyuk and J. H. Koolen, A characterization of the Grassmann graphs, J. Combin. Theory Ser. B 171 (2025), 1--27.

\bibitem{godsil} C. D. Godsil, Geometric distance-regular covers, New Zealand J. Math. 22 (1993), 31--38.

\bibitem{GR2001} C. D. Godsil and G. Royle, Algebraic Graph Theory, Springer-Verlag, Berlin, 2001.

\bibitem{-m} J. H. Koolen and S. Bang, On distance-regular graphs with smallest eigenvalue at least $-m$, J. Combin. Theory Ser. B 100 (2010), 573--584.

\bibitem{kgly24} J. H. Koolen, H. Ge, C. Lv, and Q. Yang, A Bose-Laskar-Hoffman theory for $\mu$-bounded graphs with fixed smallest eigenvalue, preprint (2024).

\bibitem{KLG2024} J. H. Koolen, C. Lv, and A. L. Gavrilyuk, A characterization of the Grassmann graphs: one missing case, preprint (2024), arXiv:2412.11018.

\bibitem{KLPY24} J. H. Koolen, C. Lv, J. Park, and Q. Yang, Bounding the intersection number $c_2$ of a distance-regular graph with classical parameters $(D, b, \alpha, \beta)$ in terms of $b$, Discrete Math. 348 (2025), no. 2, Paper No. 114239, 6 pp.

\bibitem{kp} J. H. Koolen and J. Park, Shilla distance-regular graphs, European J. Combin. 31 (2010), 2064--2073.

\bibitem{LTK2021} X. Liang, Y.-Y. Tan, and J. H. Koolen, Thin distance-regular graphs with classical parameters $(D,q,q,\frac{q^t-1}{q-1}-1)$ with $t > D$ are the Grassmann graphs, Electron. J. Combin. 28 (2021), Paper No. 4.45, 21 pp.

\bibitem{metsch2} K. Metsch, Improvement of Bruck's completion theorem, Des. Codes Cryptogr. 1 (1991), 99--116.

\bibitem{metsch} K. Metsch, On a characterization of bilinear forms graphs, European J. Combin. 20 (1999), 293--306.

\bibitem{ShYa} E. Shult and A. Yanushka, Near $n$-gons and line systems, Geom. Dedicata 9 (1980), 1--72.

\bibitem{tian2024} Y. Tian, C. Lin, B. Hou, L. Hou, and S. Gao, Further study of distance-regular graphs with classical parameters with $b < -1$, Discrete Math. 347 (2024), no. 3, Paper No. 113817, 7 pp.

\bibitem{weng} C. Weng, Classical distance-regular graphs of negative type, J. Combin. Theory Ser. B 76 (1999), 93--116.


\end{thebibliography}
\end{document}